\documentclass{amsart}

\usepackage{amssymb,amsmath,amsthm,amscd}%--pas necessaire avec amsart

\usepackage[active]{srcltx} % ajout Kristian
\usepackage[frenchb]{babel}
\usepackage{pstricks}
\usepackage[T1]{fontenc}
\usepackage{mathrsfs}

\newcommand{\C}{\mathbb{C}}%------------------------Complex numbers
\newtheorem{theo}{Theorem}[section]%
\newtheorem{prop}[theo]{Proposition}%
\newtheorem{lem}[theo]{Lemma}%

\begin{document}
\Large
%opening
\title[The Bohr Radius for an Elliptic Condenser]{The Bohr Radius for an Elliptic Condenser }
\author{Patrice LassËre \& Emmanuel Mazzilli\\ \today}

\begin{abstract}  We compute the exact value of the Bohr radius associated to an elliptic condenser of the complex plane and its Faber polynomial basis.
\end{abstract}

\keywords{Functions of a complex variable, Inequalities, Schauder basis.}
\subjclass{30B10, 30A10.}

\address{LassËre Patrice : Institut de MathÈmatiques, , UMR CNRS 5580,
Universit\'e Paul Sabatier, 
118 route de Narbonne, 31062 TOULOUSE.}
\email{lassere@math.ups-tlse.fr}

\address{Emmanuel Mazzilli : UniversitÈ Lille 1, Villeneuve d'Ascq, 59655 Cedex.}
\email{Emmanuel.Mazzilli@math.univ-lille1.fr}

\maketitle

\section{Introduction}

\bigskip

Bohr's classical theorem  \cite{bohr} asserts that if $f(z)=\sum_{n\geq 0}\,a_n z^n$ is holomorphic on the unit disc $\mathbb D$  and if $ \left \vert f(z)\right\vert<1,\ \forall\,z\in\mathbb D$ then $\sum_{n\geq 0}\, \left\vert a_n z^n\right\vert<1,\ \forall\,z\in D(0,1/3)$, 
and the constant  $1/3$ is optimal. 

\bigskip
In a previous work \cite{lasseremanu} we study the Bohr's phenomenon  in the following context  :
let    $K\subset\mathbb C$ be a continuum\footnote{i.e.  $K$ is a compact in $\mathbb C$ including at least two points, and $\overline{\mathbb C}\setminus K$ is simply connected},    
$\Phi\ :\ \overline{\mathbb C}\setminus K\to \overline{\mathbb C}\setminus{\overline {\mathbb D}}$ the unique conformal mapping satisifiying 
$\Phi(\infty)=\infty,\quad \Phi'(\infty)=\gamma>0$, and  $(F_{K,n})_n$ the sequence of its Faber polynomials (\cite{suetin}). 

\noindent This is a classical fact \cite{suetin} that  $(F_{K,n})_n$ is a Schauder basis\footnote{this means that for all $f\in\mathscr O(\Omega_{K,\rho})$ there exists an unique sequence  $(a_n)_n$ of complex numbers such that $f=\sum_{n\geq 0} a_n\varphi_n$ for the usual compact convergence of $O(\Omega_{K,\rho})$ and the same is true in $\mathscr O(K)$ equiped with its usual inductive limit topology. }   for all the spaces $ \mathscr O(\Omega_{K,\rho}),\ (\rho>1)$ and also  $ \mathscr O(K)$.  
We prove  (\cite{lasseremanu}, theorem 3.1)   that the family $(K,\Omega_{K,\rho}, (F_{K,n})_n)$ satisfies the Bohr phenomenon in the following sense : there exists $\rho_0>1$ such that  for all $\rho>\rho_0$, for all $f=\sum_n \,a_nF_{K,n}\in\mathscr O(\Omega_{K,\rho})$, if $\vert f(z)\vert<1$ for all $z\in\Omega_{K,\rho}$, then $\sum_n \,\vert a_n\vert\cdot\Vert F_{K,n}\Vert_K<1$. The infimum $\rho_K$ of all such $\rho_0$ will be called the \textbf{Bohr radius} of $K$. 

\medskip
For example, the Faber polynomial basis for the compact $\overline{ D(0,1)}$ is precisely the Taylor basis i.e. $F_{K,n}(z)=z^n$ and the levels sets are the discs $\Omega_{\overline{ D(0,1)},\rho}=D(0,\rho),\ (\rho>1)$ ; then, thanks to the classical Bohr theorem, we have a Bohr phenomenon and the Bohr radius of  $K=\overline {D(0,1)}$ is $\rho_K=3$.
 
\medskip
The particular cases $K:=[-1,1]\subset\C$ is one of the very few more examples  (see \cite{suetin}, which is the definite reference on this subject) where  the explicit form of the conformal map 
$\Phi \ :\ \Omega:=\overline{\C}\setminus K\ni z\mapsto w=\Phi(z)\in\overline{\C}\setminus \overline{D(0,1)}$
give us more  precises estimations.
In this simple case $\Phi^{-1}(w)= (w+{1/ w})/2$ is the famous Zhukovskii function.
Faber polynomials  $(F_{K,n})_n$ form a common basis of the spaces $\mathscr O(K)$ and  $\mathscr O(\Omega_{K,\rho})$, $(\rho>1)$ where the boundary $\partial\Omega_{K,\rho}$ of the level sets   are  ellipses with focus $1$ and $-1$, and excentricity
$\varepsilon={2\rho\over 1+\rho^2}$. That's why we will speak of \textbf{elliptic condenser}.

\medskip
In \cite{kap},  H.T. Kaptanoglu  \& N. Sadik study the Bohr phenomenon (with a slightly different approach) in the case of an elliptic condenser   and obtain an estimation of its Bohr radius. Their paper inspired our works and in the present one we compute the exact value of this radius. We also compute the exact value of the radius for holomorphic functions with only real coefficients in their Faber expansion. Note that  the observation that these two radius can be different (contrary to the classical Bohr's theorem) seems to be new.

\bigskip
\section{ The Sketch of the Proof and Technicals notations}

\bigskip
\noindent \textbf{The Sketch of the Proof  : } In the proof of the classical Bohr theorem, the main ingredient for an upper estimation of the Bohr radius are Carath\'eodory inequality :  

\medskip
 \textit{" let $f(z)=\sum_na_nz^n\in\mathscr O(D(0,1))$. If $\texttt{re}(f(z))>0$ for all $z\in D(0,1)$ then $\vert a_n\vert\leq 2\texttt{re}(a_0)$ for all $n\geq 1$."}
 
\medskip
In the elliptic case, the procedure is the same.  In \cite{lasseremanu} we  already prove  the following elliptic-Carath\'eodory's inequality :

 \bigskip
 \textit{ "(\cite{lasseremanu} prop 2.1)Let $f(z)=a_0+\sum_{1}^{\infty}a_nF_{K,n}(z)\in \mathscr O(\Omega_{K,\rho})$ and suppose that
$\texttt{re}({f})>0$. Then : 
$\vert a_n\vert\leq {2\texttt{re}(a_0)\over
\rho^n- \rho^{-n}},\quad\forall\,n\in\mathbb N^\star$.
Moreover,  if  $f(z)=a_0+\sum_{1}^{\infty}a_nF_{K,n}(z) \in \mathscr O(\Omega_{K,\rho})$  satisfies $\vert
f\vert <1$ with $a_0>0$, then we have\footnote{ Note that $a_0<1$ because $\vert
f\vert <1$.}  :
$\vert a_n\vert\leq {2(1-a_0)\over \rho^n-\rho^{-n}}.$ for all $n\geq 1$. "}

 %\footnote{La version prÈcÈdente Ètait ecrite dans la langue des $w$ : Let
%$f(w)=a_0+\sum_{1}^{\infty}a_n(w^n+ w^{-n})\in \mathscr O(\{1<\vert w\vert <R\})$ and suppose that
%$\texttt{re}({f})>0$. Then : 
%$$\vert a_n\vert\leq {2\texttt{re}(a_0)\over
%R^n- R^{-n}},\quad\forall\,n\in\mathbb N^\star.  $$
%Moreover,  if  $f(w)=a_0+\sum_{n=1}^{\infty}a_n(w^n+w^{-n}) \in \mathscr O(\{1<\vert w\vert <R\})$  satisfies $\vert
%f\vert <1$ with $a_0>0$, then we have  :
%$$\vert a_n\vert\leq {2(1-a_0)\over R^n-R^{-n}}.$$}

\bigskip
to deduce that  :

\bigskip
\textit{  "(\cite{lasseremanu} prop 2.3)  The elliptic condenser $(K:=[-1,1], \Omega_{K,\rho}, (F_{K,n})_n)$ satisfyies Bohr's phenomenon for all $\rho\geq \rho_0=5,1284...$." }
  
\bigskip
   which gives  an excentricity $\varepsilon_0=0.3757...$. This is already better than  $\varepsilon_0=0.373814...$($\rho_0=5,1573...$) find by   H.T. Kaptanoglu  \& N. Sadik in \cite{kap}. 

To find the exact value of $\rho_K$ we will first (in paragraph 3) prove better elliptic Caratheodory's inequalities. Then (paragraph 4) use these inequalities to get an upper bound for $\rho_K$. Finaly in paragraph 5, we prove that the upper bound obtained in the previous paragraph is optimal thanks to explicit test functions.

\bigskip
Our proof is rather technical, so before going into it, let us state clearly the mains tools we will use.

\bigskip
\noindent \textbf{Technical observations and notations  : }  First, it is  fondamental to observe that the expression of  $F_{K,n}$ is very more convenient in goal coordinates \og $w=\Phi(z)$\fg{} than in source coordinate \og $z$\fg{} because in \og$ w$\fg{} coordinates we have
$F_{K,n}(w)=w^n+ w^{-n},\  \forall\,\vert w\vert=\vert\Phi(z)\vert>1$, 
 so $\Vert F_{K,n}\Vert_{\Omega_{K,\rho}}=\rho^n+\rho^{-n},\ \forall\,\rho>1$. 

%\bigskip
%Let $\mathscr E$ denote an non degenerate ellipse with foci $\pm 1$, i.e. a level set $\Omega_R,\ (R>1) $ of the biholomorphism $\Phi \ :\ \Omega:=\overline{\C}\setminus K\ni z\mapsto w=\Phi(z)\in\{\vert w\vert >1\}$. Then, the biholomorphism $\Phi_1 \ :\ \overline{\C}\setminus \Omega_R\to\mathbb C\setminus\overline{D(0,1)}$ is of course $\Phi_1=\Phi/R$ so :
%$$F_{\mathscr E,n}(w)= F_{K,n}(w/R)=R^{-n}(w^n+ R^{2n}w^{-n}).$$
%And $\Phi_1$ extends across $\partial\mathscr E$ up to $\mathbb C\setminus K$ : $K$  is relatively to $\mathscr E$ the level set $1/R$ for $\Phi_1$ (or more precisely $K=\bigcap_{1/R<\alpha} \Omega_\alpha$). I.e. the condenser $(\mathscr E,\Phi_1)$ as level sets $\Omega_\alpha$ with $1/R<\alpha<+\infty$ and $\Omega_1=\mathscr E$.
%We have for all $1/R<\alpha$ and $\forall\,f\in\mathscr O(\Omega_{\alpha,\mathscr E})$
%$$f(z)=\sum_{n\geq 0} a_nF_{K,n}(z)=\sum_{n\geq 0} a_nR^{-n}F_{\mathscr E,n}(z)=\sum_{n\geq 0} a_n\cdot R^{-n}(w^n+ R^{2n}w^{-n}),$$
%$\forall\,z\in\Omega_{\alpha,\mathscr E}$ or $w\in 1/R<\vert w\vert<\alpha$.

\bigskip
To be in the same spirit\footnote{Of course we could have done the same by looking for $\rho>1$ such that  all  $f=\sum_{n}a_n F_{n,K}\in \mathscr O(\Omega{K,\rho})$ satisfies $\sum_n \vert a_n\vert\cdot\Vert F_{n,K}\Vert_K<1$.} that the seminal's work of   H.T. Kaptanoglu  \& N. Sadik,   and compute the Bohr radius $\rho_K$ we will procede as follow :
From now,   $\mathscr E$ will denote the domain bounded by a non degenerate ellipse with foci $\pm 1$, i.e. a level set $\Omega_{K,\rho},\ (\rho>1) $ of the biholomorphism.
 $\Phi_K \ :\ \overline{\C}\setminus K\mapsto \mathbb C\setminus\overline{D(0,1)}$. 
 
 So, the biholomorphism  $\Phi_\mathscr E \ :\  \overline{\C}\setminus \mathscr E\mapsto \mathbb C\setminus\overline{D(0,1)}$ is $\Phi_\mathscr E=\Phi_K/\rho$ which extends as a biholomorphism up to $\mathbb C\setminus K$. In another words, the level sets $\Omega_{\mathscr E,r}$ of $\Phi_\mathscr E$ are defined not only for $1<r<+\infty$ but for $1/\rho<r<+\infty$. And we have ( \cite{kap}) with $R:=\rho^{-1}$ : 
 $$F_{\mathscr E,n}(w)= w^n+ \rho^{-2n}w^{-n}=w^n+R^{2n}/w^n,\quad  R=\rho^{-1}<\vert w\vert .$$
So, for all $f \in\mathscr O(\Omega_{\mathscr E,r})$ we will have
$$f(w)=\sum_{n\geq 0}a_n\cdot F_{\mathscr E,n}(w)=\sum_{n\geq 0}a_n\cdot(w^n+\rho^{-2n}w^{-n}),\quad  r>\vert w\vert>\rho^{-1}=R.$$
 Then, following H.T. Kaptanoglu  \& N. Sadik, we are going to look for the largest $0<R<1$ such that we have a bohr phenomenon for the family $(K, \mathscr E=\Omega_{K,R^{-1}}, (F_{K,n})_n$. If we note $R_{B}$ this largest $R$, clearly $\rho_K={1/ R_{B}}$.

\vfill\eject
\bigskip
\section{\og Elliptic\fg{} Caratheodory Inequalities. } 

\bigskip
%Let  $\mathscr E$ denote an non degenerate ellipse with foci $\pm 1$, i.e. a level set $\Omega_R,\ (R>1) $ of the biholomorphism. Then, 

%--  Up to a rotation, we always suppose  in this paragraph  that  $a_0>0$.

%-- Remember  that the bihomorphism  allow us to translate all computations on the ellipse onto the discs via the formula
%$$\forall\,1>\vert w\vert>R\ :\quad F_{\mathscr E,n}(\Phi^{-1}(w))=w^n+R^{2n}w^{-n}.$$

Let $f=\sum_{n\geq 0} a_n  F_{\mathscr E,n}\in\mathscr O(\mathscr E)$, up to a rotation, we always suppose  in this paragraph  that  $a_0>0$. Then elementary computations gives for all $n\geq 1$  :

\begin{small}
 \begin{align}
 \int_0^{2\pi} f\left(\Phi^{-1}(e^{i\theta})\right)e^{in\theta}d\theta&=\int_0^{2\pi}\left( a_0+ \sum_{k\geq 1} a_k
\left( e^{ik\theta}+R^{2k}e^{-ik\theta}\right)\right) e^{in\theta}d\theta=R^{2n}a_n\\
  \int_0^{2\pi} \overline{f(\Phi^{-1}(e^{i\theta}))}e^{in\theta}d\theta&=\int_0^{2\pi}\left( \overline{a_0}+ \sum_{k\geq 1} \overline{a_k}
\left( e^{-ik\theta}+R^{2k}e^{ik\theta}\right)\right) e^{in\theta}d\theta=\overline{a_n},
\end{align} 
\end{small}

specially :
\begin{small}
 \begin{align}
  R^{4n}a_{2n}=\int_0^{2\pi} f\left(\Phi^{-1}(e^{i\theta})\right)e^{2in\theta}d\theta,\quad 
a_{2n}=           \int_0^{2\pi} f\left(\Phi^{-1}(e^{i\theta})\right)e^{-2in\theta}d\theta,\\
  R^{4n}\overline{a_{2n}}=\int_0^{2\pi} \overline{f\left(\Phi^{-1}(e^{i\theta})\right)}e^{-2in\theta}d\theta,\quad
\overline{a_{2n}}= \int_0^{2\pi} \overline{f\left(\Phi^{-1}(e^{i\theta})\right)}e^{2in\theta}d\theta.
 \end{align}
\end{small}

Our goal in this paragraph is two prove the following "elliptic Caratheodory's type  inequality" : 

\bigskip
\begin{prop} Let $f=\sum_{n\geq 0} a_n  F_{\mathscr E,n}\in\mathscr O(\mathscr E)$. If \texttt{re}$(f)\geq 0$ on $\mathscr E$ and $R\leq 0.2053...$, then for all $n\in\mathbb N^\star$   : 
$$ \vert a_{n}\vert R^n+\vert a_{2n}\vert R^{2n} 
\leq \dfrac{2\texttt{re}(a_0)R^n}{1-R^{2n}}+ \dfrac{2\texttt{re}(a_0)R^{2n}}{1+R^{4n}}.
 $$

\end{prop}

\bigskip
First we need  two lemmas :

\bigskip
\begin{lem} \textbf{("Classical"  Caratheodory's inequality)} 

Let $f(z)=\sum_{n\geq 0} a_n  F_{\mathscr E,n}(z)\in\mathscr O(\mathscr E)$. If \texttt{re}$(f)\geq 0$ on $\mathscr E$, then for all $n\in\mathbb N^\star$ : 
\begin{equation}\left\vert R^{2n}a_n+\overline{a_n}\right\vert=(1+R^{2n})^2 \texttt{re}^2(a_n)+(1-R^{2n})^2 \texttt{im}^2(a_n) \leq 4\texttt{re}^2(a_0).\end{equation}
Particulary :
\begin{equation}\left\vert \texttt{re}(a_n)\right\vert \leq \dfrac{2\texttt{re}(a_0)}{1+R^{2n}},\quad\forall\,n\in\mathbb N^\star.\end{equation}
\end{lem}

\bigskip
\noindent\textbf{Proof : } For all $n\in\mathbb N^\star$ with  (1) and (2) we have : 
\begin{small}
 $$ 
 \left\vert R^{2n}a_n+\overline{a_n}\right\vert = \left\vert \int_0^{2\pi} 2\texttt{re}\left( f\left(\Phi^{-1}(e^{i\theta})\right)\right) e^{in\theta}d\theta \right\vert 
 \leq 2 \int_0^{2\pi} 2\texttt{re}\left( f\left(\Phi^{-1}(e^{i\theta})\right)\right) d\theta 
 =2\texttt{re}(a_0). 
 $$
\end{small}
 
\noindent So  for all $n\in\mathbb N^\star$ : 
$$\left\vert R^{2n}a_n+\overline{a_n}\right\vert=(1+R^{2n})^2 \texttt{re}^2(a_n)+(1-R^{2n})^2 \texttt{im}^2(a_n)\leq 4 \texttt{re}^2(a_0),$$
i.e.
$$\vert\texttt{re}(a_n)\vert\leq \dfrac{\sqrt{4 \texttt{re}^2(a_0)-(1-R^{2n})^2 \texttt{im}^2(a_n)}}{1+R^{2n}} \leq \dfrac{2\texttt{re}(a_0)}{1+R^{2n}}.$$
QED.\footnote{Observe that, if all the $a_n\in\mathbb R$ then we have the stronger inequalities $\vert a_n\vert\leq \frac{2 \texttt{re}(a_0)}{1+R^{2n}},\ n\geq 1$ which will be fundamental for the "real" case.}\hfill$\blacksquare$

\bigskip
\begin{lem} Under the lates assumptions we have for all $n\in\mathbb N^\star$ :
$$\vert a_n\vert \leq 2\cdot\dfrac{\sqrt{1+R^{4n}}}{1-R^{4n}}\sqrt{ \texttt{re}(a_0)-R^{2n}\texttt{re}(a_{2n}) }\cdot \sqrt{\texttt{re}(a_0)}.$$
\end{lem}

\bigskip
\noindent\textbf{Proof  : } It will be more convenient in the sequel to state  
$$\Theta_n(e^{i\theta}):= e^{in\theta}-R^{2n}e^{-in\theta}.$$
Then   $\vert \Theta_n(e^{i\theta})\vert = 1+R^{4n}-R^{2n}\left( e^{-2in\theta}+e^{2in\theta}\right)$. 
So, with (3) and (4), we can write :
{%\small 
 $$\begin{aligned}
2\int_0^{2\pi} \texttt{re}\left(f\left(\Phi^{-1}(e^{i\theta})\right)\right)&\cdot\vert   \Theta_n(e^{i\theta})\vert d\theta  =2(1+R^{4n})\int_0^{2\pi} \texttt{re}\left(f\left(\Phi^{-1}(e^{i\theta})\right)\right)d\theta\\
&  -R^{2n}\int_0^{2\pi}\,2\texttt{re}\left(f\left(\Phi^{-1}(e^{i\theta})\right)\right)\left( e^{-2in\theta}+e^{2in\theta}\right)d\theta\\
&=2(1+R^{4n})\texttt{re}(a_0)-R^{2n}\left[ a_{2n}+\overline{a_{2n}}+R^{4n}(a_{2n}+\overline{a_{2n}})\right]\\
&=2(1+R^{4n})\left[ \texttt{re}(a_0)-R^{2n}\texttt{re}(a_{2n})\right].
\end{aligned}$$}
And also :
{%\small
 $$\begin{aligned}\int_0^{2\pi}\,f\left(\Phi^{-1}(e^{i\theta})\right)\cdot \Theta_n(e^{i\theta})d\theta
%&=\int_0^{2\pi} \left( a_0+\sum_{k\geq 1} a_k \left(e^{ik\theta}+R^{2k}e^{-ik\theta}\right)\right)\cdot \left(e^{in\theta}-R^{2k}e^{-in\theta}\right)d\theta\\
& =0,\\
\int_0^{2\pi}\,\overline{f\left(\Phi^{-1}(e^{i\theta})\right)} \cdot\Theta_n(e^{i\theta})d\theta
%&=\int_0^{2\pi} \left( \overline{a_0}+\sum_{k\geq 1} \overline{a_k} \left(e^{-ik\theta}+R^{2k}e^{ik\theta}\right)\right)\cdot \left(e^{in\theta}-R^{2k}e^{-in\theta}\right)d\theta \\
&= \overline{a_n}(1-R^{4n}).
\end{aligned}$$}
Consequently
 {\small $$\begin{aligned}
& \left\vert \int_0^{2\pi}\, \left( \overline{f\left(\Phi^{-1}(e^{i\theta})\right)}+f\left(\Phi^{-1}(e^{i\theta})\right) \right)  \Theta_n(e^{i\theta}) d\theta \right\vert \\
&\hspace{1cm} \leq 2\int_0^{2\pi}  \texttt{re}\left(f\left(\Phi^{-1}(e^{i\theta})\right)\right)\vert\Theta_n(e^{i\theta})\vert d\theta\\
%& =
%2\int_0^{2\pi} \sqrt{ \texttt{re}\left(f\left(\Phi^{-1}(e^{i\theta})\right)\right)}\cdot\vert\Theta_n(e^{i\theta})\vert \cdot \sqrt{ \texttt{re}\left(f\left(\Phi^{-1}
%(e^{i\theta})\right)\right)} d\theta\\
 &\hspace{1cm}  \leq
2\left( \int_0^{2\pi} \texttt{re}\left(f\left(\Phi^{-1}(e^{i\theta})\right)\right)\cdot\vert\Theta_n(e^{i\theta})\vert^2 d\theta \right)^{1/2}\cdot
\left( \int_0^{2\pi}  \texttt{re}\left(f\left(\Phi^{-1}(e^{i\theta})\right)\right) d\theta \right)^{1/2}
\end{aligned}$$}
thanks to Cauchy-Schwarz. With the three previous identities  we gets :
$$\vert {a_n}\vert\cdot (1-R^{4n})\leq 2\sqrt{\texttt{re}(a_0)\cdot(1+R^{4n})\cdot\left( \texttt{re}(a_0)-R^{2n}\texttt{re}(a_{2n})\right)}$$
for all $n\geq 1$,  the desired inequality. QED.\hfill$\blacksquare$

\bigskip
Now we are able to give the 

\bigskip
\noindent\textbf{Proof of the proposition 3.1 : } The "classical Caratheodory inequality" (lemma 3.1)  gives  for  $a_{2n}$
$$(1+R^{4n})^2 \texttt{re}^2(a_{2n})+(1-R^{4n})^2 \texttt{im}^2(a_{2n}) \leq 4\texttt{re}^2(a_0),$$
which implies
$$\texttt{im}^2(a_{2n})\leq \dfrac{4\texttt{re}^2(a_0)-(1+R^{4n})^2 \texttt{re}^2(a_{2n})}{(1-R^{4n})^2}$$
  so
$$\begin{aligned}
\vert a_{2n}\vert^2 &\leq 
\dfrac{4\texttt{re}^2(a_0)}{(1-R^{4n})^2}+\texttt{re}^2(a_{2n})\left[1-\dfrac{(1+R^{4n})^2}{(1-R^{4n})^2}\right]\\
&\leq \dfrac{4}{(1-R^{4n})^2}\cdot \left[ \texttt{re}^2(a_0)-R^{4n}\texttt{re}^2(a_{2n})\right]
\end{aligned}$$
or :
\begin{equation}\vert a_{2n}\vert\leq \dfrac{2}{1-R^{4n}}\cdot \sqrt{\texttt{re}^2(a_0)-R^{4n}\texttt{re}^2(a_{2n})},\end{equation}
for all $n\in\mathbb N^\star$. This last inequality associated with lemme 3.2 lead us to the main estimation :
{\small $$\begin{aligned}\vert a_{n}\vert R^n&+\vert a_{2n}\vert R^{2n} \leq \\
&\leq \dfrac{2R^n}{1-R^{4n}}\left[ 
\sqrt{\texttt{re}(a_0)(1+R^{4n})(\texttt{re}(a_0)-R^{2n}\texttt{re}(a_{2n}))}+R^n\sqrt{\texttt{re}^2(a_0)-R^{4n}\texttt{re}^2(a_{2n})}\right] \\
&\leq 
\dfrac{2R^n}{1-R^{4n}}\left[ 
\sqrt{\texttt{re}(a_0)(1+R^{4n})(\texttt{re}(a_0)+R^{2n}\texttt{re}(a_{2n}))}+R^n\sqrt{\texttt{re}^2(a_0)-R^{4n}\texttt{re}^2(a_{2n})}\right]\\
&=\dfrac{2R^n}{1-R^{4n}}\left[ 
\sqrt{\texttt{re}(a_0)(1+R^{4n})(\texttt{re}(a_0)+R^{2n}x)}+R^n\sqrt{\texttt{re}^2(a_0)-R^{4n}x^2}\right]\\
&:=\dfrac{2R^n}{1-R^{4n}}G(x)
\end{aligned}$$}
where $x=\vert \texttt{re}(a_{2n})\vert\in [0,\frac{2\texttt{re}(a_{0})}{1+R^{4n}}]$ . Now, let us  maximize $G(x)$ on $[0,\frac{2\texttt{re}(a_{0})}{1+R^{4n}}]$ : 

$$\begin{aligned}
G(x) &= \sqrt{\texttt{re}(a_0)(1+R^{4n})(\texttt{re}(a_0)+R^{2n}x)}+R^n\sqrt{\texttt{re}^2(a_0)-R^{4n}x^2} \\
G'(x)&= \dfrac{R^{2n}}{2\sqrt{\texttt{re}(a_0)+R^{2n}x}}\left( \sqrt{\texttt{re}(a_0)(1+R^{4n})}-\dfrac{2R^{3n}x}{\sqrt{\texttt{re}(a_0)-R^{2n}x}}\right).
\end{aligned}$$
So
$$%\begin{aligned}
G'(x)=0\ \ 
%&\iff\ \ 2R^{3n}x = \sqrt{\texttt{re}(a_0)(1+R^{4n})}\cdot \sqrt{\texttt{re}(a_0)-R^{2n}x}\\
%&\iff\ \ 4R^{6n}x^2 = \texttt{re}(a_0)(1+R^{4n})\cdot (\texttt{re}(a_0)-R^{2n}x)\\
%&\iff\ \ 4R^{6n}x^2 = \texttt{re}(a_0)\left( \texttt{re}(a_0) -R^{2n}x-R^{6n}x+\texttt{re}(a_0)R^{4n}  \right)\\
%&\iff\ \ 4R^{6n}x^2 = \texttt{re}^2(a_0)(1+R^{4n})-\texttt{re}(a_0)R^{2n}x(1+R^{4n})\\
\iff\ \ 4R^{6n}x^2 +\texttt{re}(a_0)R^{2n}(1+R^{4n})x-\texttt{re}^2(a_0)(1+R^{4n})=0
%\end{aligned}
$$
whose roots are
$$\begin{cases}
x_1&=-\dfrac{\texttt{re}(a_0)\sqrt{1+R^{4n}}}{8R^{4n}}\left[\sqrt{1+R^{4n}+16R^{2n}}+\sqrt{1+R^{4n}}\right]<0,\\
x_2&=\dfrac{\texttt{re}(a_0)\sqrt{1+R^{4n}}}{8R^{4n}}\left[\sqrt{1+R^{4n}+16R^{2n}}-\sqrt{1+R^{4n}}\right]>0.
\end{cases}$$
Because $G'(x)\geq 0$ for $x\in [0,x_2]$ and $G'(x)<0$ for $x>x_2$, we have to study the sign of  $x_2-\frac{2\texttt{re}(a_0)}{1+R^{4n}}$. First, observe that\footnote{because
$\frac{1}{2\sqrt{a}}(b-a)\geq \sqrt{b}-\sqrt{a}=\int_a^b\,\frac{dt}{2\sqrt{t}}\geq \frac{1}{2\sqrt{b}}(b-a),\quad \forall\, a<b\in\mathbb R_+^\star$.}
$$\sqrt{1+R^{4n}+16R^{2n}}-\sqrt{1+R^{4n}}\geq \dfrac{1}{2}\cdot\dfrac{16R^{2n}}{\sqrt{1+R^{4n}+16R^{2n}}},$$
So
\begin{small}
$$\begin{aligned}
x_2-\frac{2\texttt{re}(a_0)}{1+R^{4n}} &= \dfrac{\texttt{re}(a_0)(\sqrt{1+R^{4n}}}{8R^{4n}}\left[\sqrt{1+R^{4n}+16R^{2n}}-\sqrt{1+R^{4n}}\right]-\frac{2\texttt{re}(a_0)}{1+R^{4n}}\\
&\geq
\dfrac{\texttt{re}(a_0)\sqrt{1+R^{4n}}}{8R^{4n}}\cdot \dfrac{1}{2}\cdot\dfrac{16R^{2n}}{\sqrt{1+R^{4n}+16R^{2n}}} -\frac{2\texttt{re}(a_0)}{1+R^{4n}}\\
&\geq
\dfrac{\texttt{re}(a_0)\sqrt{1+R^{4n}}}{R^{2n}\sqrt{1+R^{4n}+16R^{2n}}} -\frac{2\texttt{re}(a_0)}{1+R^{4n}}\\
&\geq 
\dfrac{\texttt{re}(a_0)\left[(1+R^{4n})\sqrt{1+R^{4n}}-2R^{2n}\sqrt{1+R^{4n}+16R^{2n}}\right]}{R^{2n}(1+R^{4n})\sqrt{1+R^{4n}+16R^{2n}}} 
\end{aligned}$$
\end{small}
Now let us study the sign of the numerator $(1+R^{4n})^{3/2}-2R^{2n}\sqrt{1+R^{4n}+16R^{2n}}$ : as we saw it just few lines above the inequality $\frac{1}{2\sqrt{a}}(b-a)\geq \sqrt{a}-\sqrt{b}$ gives
$$\sqrt{1+R^{4n}+16R^{2n}}\leq \sqrt{1+R^{4n}}+\dfrac{8R^{2n}}{\sqrt{1+R^{4n}}},$$
which implies
$$\begin{aligned}
(1+R^{4n})^{3/2}-2R^{2n}&\sqrt{1+R^{4n}+16R^{2n}}\geq \\
&\geq (1+R^{4n})^{3/2}-2R^{2n} \sqrt{1+R^{4n}}-\dfrac{16R^{4n}}{\sqrt{1+R^{4n}}}\\
&\geq 
\sqrt{1+R^{4n}} (1+R^{4n}-2R^{2n})-\dfrac{16R^{4n}}{\sqrt{1+R^{4n}}}\\
&\geq
\sqrt{1+R^{4n}} (1-R^{2n})^2-{16R^{4n}}\\
&\geq (1-R^{2n})^2-{16R^{4n}}=-15R^{4n}-2R^{2n}+1.
\end{aligned}$$
But, $-15R^{4n}-2R^{2n}+1\geq 0$ if $R^{2n}\in [-1/3,1/5]$ which is the case if $R^2\leq 1/5$ i.e. $0\leq R\leq 1/\sqrt{5}\simeq 0.447...$ which is more than confortable because\footnote{Remember that if  $R\geq 0,2053...$ we have no Bohr's phenomenon as we saw it in the last paragraph.} $R\leq 0.2053...$  

\bigskip
So, for $R\leq 0.2053...$ we have for all $n\in\mathbb N^\star$ : 
$$\begin{aligned}\vert a_{n}\vert R^n+\vert a_{2n}\vert R^{2n} 
&\leq \dfrac{2R^n}{1-R^{4n}}G\left( \dfrac{2\texttt{re}(a_0)}{1+R^{4n}}  \right)\\
&\leq \dfrac{2\texttt{re}(a_0)R^n}{1-R^{2n}}+ \dfrac{2\texttt{re}(a_0)R^{2n}}{1+R^{4n}}
\end{aligned}$$
Which is better that the expected estimation :
$$\vert a_n\vert\leq \dfrac{2\texttt{re}(a_0)}{1-R^{2n}},\quad \forall\,n\in\mathbb N^\star,\ a_n\not\in\mathbb R.$$
\hfill$\blacksquare$

 \vfill\eject
\section{Minoration for the Bohr radius}

\bigskip

Remember the notations : $\Theta_n(e^ {i\theta})=e^{in\theta}-R^{2n}e^{-in\theta}$,
$F_{\mathscr E,n}(\Phi^{-1}(e^{i\theta}))=e^{in\theta}+R^{2n}e^{-in\theta}$ is $n$-th Faber's polynomial for the ellipse $\mathscr E$. Let $f=\sum a_nF_{\mathscr E,n}\in\mathscr O(\mathscr E)$ with (without loosing any generality) a positive real part and $a_0>0$. Then we have for all $n\geq 1$  :  
%{\color{blue}  (ici on intËgre sur le bord de l'ellipse $\mathscr E$ ce qui n'est pas lÈgal ‡ priori, ce n'est pas trËs grave, mais il faut en dire un mot : par exemple, expliquer qu'‡ priori on devrait faire $f\in O( \Omega_{\mathscr E,r})$ avec $r=1+\varepsilon$ puis epsilon tends vers $0$...)}
$$\begin{aligned}
\int_{0}^{2\pi} f\left(\Phi^{-1}(e^{i\theta})\right)\Theta_n^{2}(e^ {i\theta})d\theta&=-2a_0R^{2n}+2a_{2n}R^{4n},\\
\int_{0}^{2\pi}\overline{ f\left(\Phi^{-1}(e^{i\theta})\right)}\Theta_n^{2}(e^ {i\theta})d\theta&=-2a_0R^{2n}+\overline{ a_{2n}}(1+R^{8n}).
\end{aligned}$$
So :
$$\int_{0}^{2\pi}(f\left(\Phi^{-1}(e^{i\theta})\right)+\overline{ f\left(\Phi^{-1}(e^{i\theta})\right)})\Theta_n^{2}(e^ {i\theta})d\theta=-4a_0R^{2n}+2a_{2n}R^{4n}+\overline{ a_{2n}}(1+R^{8n})$$
which implies :
 $$\left\vert
\int_{0}^{2\pi}\left(f\left(\Phi^{-1}(e^{i\theta})\right)+\overline{ f\left(\Phi^{-1}(e^{i\theta})\right)}\right)\Theta_n^{2}(e^ {i\theta})d\theta\right\vert^2\geq \left(-4a_0R^{2n}+\texttt{re} (a_{2n})(1+R^{4n})^2\right)^2.$$
On the other side, using  Cauchy-Schwarz's inequality (remember that $\texttt{re}(f)\geq 0$), we have
$$\begin{aligned}
&\left\vert
\int_{0}^{2\pi}(f\left(\Phi^{-1}(e^{i\theta})\right)+\overline{ f\left(\Phi^{-1}(e^{i\theta})\right)})\Theta_n^{2}(e^ {i\theta})d\theta\right\vert \\
&\qquad\qquad \leq \int_{0}^{2\pi}\left\vert f\left(\Phi^{-1}(e^{i\theta})\right)+\overline{ f\left(\Phi^{-1}(e^{i\theta})\right)}\right\vert\cdot\vert\Theta_n(e^ {i\theta})\vert^2 d\theta\\
&\qquad\qquad \leq\left(
\int_{0}^{2\pi}2\texttt{re} \left(f\left(\Phi^{-1}(e^{i\theta})\right)\right)\vert\Theta_n(e^ {i\theta})\vert^4d\theta\right)^{{1/2}} 
\left(\int_ {0}^{2\pi}2 \texttt{re}\left(f\left(\Phi^{-1}(e^{i\theta})\right)\right)d\theta\right)^{{1/ 2}}.
\end{aligned}$$
Easy computation gives 
$$\vert
\Theta_n(e^ {i\theta})\vert^4=(1+4R^{4n}+R^{8n})-2R^{2n}(1+R^{4n})(e^{2in\theta}+e^{-2in\theta})+R^{4n}(e^{4in\theta}+e^{-4in\theta}),$$
so
$$\begin{aligned}&\int_{0}^{2\pi}2\texttt{re}\left(f(\Phi^{-1}(e^{i\theta})\right)\vert\Theta_n(e^ {i\theta})\vert^4d\theta\\
&\qquad=
2(1+R^{4n}+R^{8n})a_0-4R^{2n}(1+R^{4n})^2\texttt{re} (a_{2n})+2R^{4n}(1+R^{8n})\texttt{re} (a_{4n}).
\end{aligned}$$
Then, we can deduce the main inequality %(thanks to Cauchy-Schwarz's one above) :
$$\texttt{re}^2(a_{2n})\leq {4a_0(1+R^{8n})\over (1+R^{4n})^4}\left(a_0+R^{4n}\texttt{re} (a_{4n})\right).$$
This implies%\footnote{plus tard les changements : now put $x:=\texttt{re}(a_{2n})$ and $y:=\texttt{re}(a_{4n})$.} 
$$ \texttt{re}(a_{4n}) \geq {1\over R^{4n}}\left({\texttt{re}(a_{2n})^{2}(1+R^{4n})^{4}\over
  4a_0(1+R^{8n})}-a_0\right):=h_n(\texttt{re}(a_{2n})).$$
\noindent Using Carath\'eodory's inequality  (7) \footnote{This is inequality (7) in the proof of proposition 3.1 page 6.} we have : $$\vert a_{4n}\vert\leq {2\over
  (1-R^{8n})}\sqrt{a_{0}^2-R^{8n}\texttt{re}(a_{4n})^2}:=g_n(\texttt{re}(a_{4n})),$$
Because of these two inequalities, define
$$\begin{aligned}
g_n(u)&={2\over  (1-R^{8n})}\sqrt{a_{0}^2-R^{8n}u^2},\quad &0\leq u=\texttt{re}(a_{4n})\leq {2a_0\over 1+R^{8n}}\\
h_n(v)&={1\over R^{4n}}\left({v^2(1+R^{4n})^ 4\over 4a_0(1+R^{8n})}-a_0
\right),\quad &0\leq v=\texttt{re}(a_{2n})\leq {2a_0\over 1+R^{4n}}.
\end{aligned}  $$
Elementary computation assures that
 $ g_n$  decrease on $[0,2a_0/(1+R^{8n})]$, $h_n$ increase on $ [0,2a_0/(1+R^{4n})]$.
with
$$\begin{aligned}
g_n\left({2a_0\over 1+R^{8n}}\right)&={2a_0\over 1+R^{8n}},\quad g'_n\left({2a_0\over 1+R^{8n}}\right)=-{4R^ {8n}\over (1-R^{8n})^ 2}\\
h_n\left({2a_0\over 1+R^{4n}}\right)&={2a_0\over 1+R^{8n}},\quad h'_n\left({2a_0\over 1+R^{4n}}\right)= {(1+R^{4n})^{3}\over R^ {4N}(1+R^{8n})}
\end{aligned}$$
And (remember that we have already $R<0.2053..$)
$$h'_n\left({2a_0\over 1+R^{4n}}\right) \leq {2\over R^{4n}},\quad \forall\,0\leq R\leq
{1\over 2}.$$

From now on, to simplify,  we will note $x_{0}^{n}:={2a_0\over 1+R^{4n}}$,.

\bigskip
\begin{lem}  Let $x_{1}$ be the unique value in $[0,x_{0}^{n}]$ such that $ h_n(x_{1})=0$. Define $\phi_n$
 on $[0,x_{0}^{n}]$  by 
$$\phi_n(t)=\begin{cases} 
 {2a_0\over 1-R^{8n}},\quad &{\text{ if }}  t\in [0, x_{1}],\\
g_n\circ h_n(t),\quad &{\text{ if }} 
 t\in [x_{1},x_{0}^{n}].
 \end{cases}$$
  Then,  $\phi_n\in C^1([0,x_{0}^{n}])$ and we have the following estimation  :
$$R^{4n}\phi'_n(t)\geq -4{R^{8n}(1+R^{4n})\over
(1-R^{4n})^{2}(1+R^{8n})},\quad\forall\,t\in[x_{1},x_{0}^{n}].$$
Moreover, if $R\leq {1\over 2}$ then : \quad  $R^{4n}\phi'_n(t)\geq -8R^{8n}$.
\end{lem}

\bigskip
\noindent\textbf{Proof : } First part is trivial. For the last one, we have
$R^{4n}(g_n\circ h_n)'(t) =R^{4n}h_{n}^{'}(t)g_{n}^{'}(h_n(t))$ and we now that :

\medskip
-- $h'_{n}$ increase on $[0,{2a_0\over 1+R^{4n}}]$ and take positive values. 

-- $g'_{n}$ decrease

-- $h_n\geq 0$  on  $x\geq x_{1}$
 
\noindent so,  $g_{n}^{'}(h_n(x))$ is negative and reach its minimum at ${2a_0\over 1+R^{4n}}$. The,  easy computation gives the required inequality.\hfill$\blacksquare$

\bigskip

We will also need the following estimations :

\begin{lem}Fix $n_0\geq 1$, then :
\begin{enumerate}
\item Fix $k\in\mathbb N$ and let $x_1:=h^{-1}_{n_0}\circ\cdots \circ h^{-1}_{2^kn_0}(0)$. Then the function defined on $[0,x_{0}^{n}]$  by 
$$\phi_{2^{k}n_0}(t)=\begin{cases} 
 {2a_0\over 1-R^{n_02^{k+3}}},\quad &{\text{ if }}  t\in [0, x_{1}],\\
g_{2^{k}n_0}\circ h_{2^kn_0}\circ h_{2^{k-1}}\circ \cdots
\circ h_{n_0}(t),\quad &{\text{ if }} 
 t\in [x_{1},x_{0}^{n}].
 \end{cases}$$
satisfies 
$$\forall\,t\in [0,x_{0}^{n_0}], R\in]0,1/2] \ :\quad  R^{2^{k+2}n_0}\phi'_{2^{k}n_0}(x)\geq -2^{k+3}R^{n_02^{k+3}}.$$  
\item $\displaystyle -\sum_{k=0}^{\infty}2^{k+3}R^{n_02^{k+3}}\geq -{8R^{8n_0}\over
  1-2R^{8n_0/3}}\geq -16R^{8n_0},\quad\forall\, R\leq 1/2.$
\end{enumerate} 
\end{lem}

\noindent\textbf{Proof : } 1) We have :
$$R^{2^{k+2}n_0}(\phi_{2^{k}n_0}(x))^{'}=R^{2^{k+2}n_0}g^{'}_{2^{k}n_0}(x_0^{2^{k}n_0})\cdot h^{'}_{2^{k}n_0}(x_0^{2^{k}n_0})\cdot h^{'}_{2^{k-1}n_0}(x_0^{2^{k-1}n_0})
\cdots h^{'}_{n_0}(x_{0}^{n_0}).$$
By the lemma 4.1 and the remarks before, we have the minoration
$$R^{2^{k+2}n_0}(\phi_{2^{k}n_0}(x))^{'}\geq -8R^{n_02^{k+3}}\times{2\over R^{n_02^{k+1}}}\times\cdots\times{2\over
  R^{n_02^{2}}}.$$
that is $R^{2^{k+2}n_0}(\phi_{2^{k}n_0}(x))^{'}\geq -2^{k+3}R^{n_02^{k+3}}$.

\medskip
\noindent 2) $n_0\geq 1$ being fixed  
$$-\sum_{k=0}^{\infty}2^{k+3}R^{n_02^{k+3}}\geq
-\sum_{k=0}^{\infty}\left(2^{(k+3)2^{-k-3}/ n_0}R\right)^{n_02^{k+3}}\geq
-\sum_{k=0}^{\infty}\left(2^{{8n_0/3}}R\right)^{n_02^{k+3}}.$$
Now,   $R\leq 1/2$ implies $2^{{3/ 8n_0}}R\leq 1$ so
$$-\sum_{k=0}^{\infty}\left(2^{{3/ 8n_0}}R\right)^{n_02^{k+3}}\geq
-\sum_{k=0}^{\infty}\left(2R^{{8n_0/ 3}}\right)^{k+3}\geq -{8R^{8n_0}\over
  1-2R^{8n_0/ 3}},$$
because $n_02^{k+3}\geq 8(k+3)n_0/ 3$. Finaly
$$-\sum_{k=0}^{\infty}2^{k+3}R^{n_02^{k+3}}\geq -{8R^{8n_0}\over
  1-2R^{8n_0/ 3}}\geq -16R^{8n_0},\quad\forall\, R\leq 1/2.$$
\hfill$\blacksquare$  

\bigskip

Because one more time of inequality (7)  we have   $n\geq 1$ : 
$$ \vert a_n\vert
R^{n}+\vert a_{2n}\vert R^{2n}\leq f_n(x),\quad\forall_,n\geq 1.$$
where 
$$f_n(x):={2R^{n}\over
  1-R^{4n}}\left(  \sqrt{a_0(1+R^{4n})(a_0-R^{2n}x)}+R^n\sqrt{a_0^2-R^{4n}x^2}\right).$$
  
\begin{lem} For all $R\leq {1\over 2}$ and $n\geq 1$, we have   :
$$\forall\,x\in\left[0,{2a_0\over 1+R^{4n}}\right]\quad :\quad f_{n}^{'}(x)\geq {R^{3n}\over 4}.$$
\end{lem}

\bigskip

\noindent\textbf{Proof : } Write $f_n(x)={R^{n}\over 1-R^{4n}}\theta_n(x).$ Then, we have
$$\theta_n^{'}(x)={1\over
  a_0^2-R^{4n}x^{2}}\bigg(\sqrt{a_0(1+R^{4n})(a_0-R^{2n}x)}-2R^{3n}x\bigg);$$
and after some computations :
$$\begin{aligned}\theta_n^{''}(x)&=-{R^{2n}\over (a_0^{2}-R^{4n}x^{2})^{2}
}\bigg(2R^{2n}x\sqrt{a_0(1+R^{4n})(a_0-R^{2n}x)}-4R^{5n}x^2\\
&+{1\over 2}\sqrt{a_0(1+R^{4n})(a_0-R^{2n}x)}(a_0+R^{2n}x)+2R^{n}(a_0^{2}-R^{4n}x^{2})\bigg)\\
&=-{R^{2n}\over (a_0^{2}-R^{4n}x^{2})^{2}}\bigg({1\over
  2}\sqrt{a_0(1+R^{4n})(a_0-R^{2n}x)}(5R^{2n}x+a_0)\\
  &\hspace{7cm} +2R^{n}(a_0^{2}-3R^{4n}x^2)\bigg).
  \end{aligned}$$
If $x_0^n={2a_0\over 1+R^{4n}}$, we have $a_0^{2}-3R^{4n}x^2={a_0^2\over
  (1+R^{4n})^2}(1-10R^{4n}+R^{8n})\geq a_0^{2}(1-10R^{4})$ 
(for all $n\geq
1$)  and so is positive if  $R\leq {1/ 2}$. So $\theta_n^{''}$ is negative on $[0,x_{0}^{n}]$ and the infimum of $f_n^{'}$ on $[0,x_{0}^{n}]$ is
${R^{n}\over 1-R^{4n}}\theta^{'}(x_{0}^{n})$. After some simple computations we get $\theta_n^{'}(x_{0}^{n})={R^{2n}\over
  1-R^{4n}}(1-R^{2n}-4R^{3n}+R^{4n}-R^{6n}),$ so,    $R\leq {1/ 2}$ implies 
$\theta_n^{'}(x_{0}^{n})\geq{ R^{2n}/ 4}$. Consequently  $f_n^{'}(x)\geq{
  R^{3n}/ 4}$ if  $R\leq{1/ 2}$,
on $[0,x_{0}^{n}]$. Q.E.D. \hfill$\blacksquare$

\bigskip
\begin{lem} We have the following elliptic version of Carth\'eodory's inequality  : let $f=a_0+\sum_{n=1}^{\infty}a_nF_{n,\mathscr E}$ be holomorphic on the ellipse $\mathscr E$. If  $\texttt{re}(f)>0$ then :
$$\sum_{n=1}^{\infty}R^{n}\vert a_n\vert\leq \sum_{n=0[2],n\geq 1}{2R^{n}a_0\over
  1+R^{2n}}+\sum_{n=1[2]}{2R^{n}a_0\over 1-R^{2n}}.$$
\end{lem}

\bigskip

\noindent\textbf{Proof : } Suppose $n_0$  odd and let $x=\texttt{re}(a_{2n_0})\in [0,x_{0}^{n_0}]$. Then because of the preceedings lemmas  
$$R^{n_0}\vert a_{n_0}\vert +R^{2n_0}\vert a_{2n_0}\vert +
\sum_{k=0}^{\infty}R^{2^{k+2}n_0}\vert a_{2^{k+2}n_0}\vert\leq f_{n_0}(x)+
\sum_{k=0}^{\infty}\phi_{2^{k}n_0}(x),$$
 The derivative of the function on right side is greater than ${R^{3n_0}\over 4}-16R^{8n_0}$ on $[0,x_{0}^{n_0}]$ ; so it is positive if  $R\leq
    0.4$. This implies that the right side of the inequality is an increasing function and so :
$$R^{n_0}\vert a_{n_0}\vert +R^{2n_0}\vert a_{2n_0}\vert+
\sum_{k=0}^{\infty}R^{2^{k+2}n_0}\vert a_{2^{k+2}n_0}\vert\leq {2R^{n_0}a_0\over
  1-R^{2n_0}}+\sum_{k=1}^{\infty}{2R^{2^{k}n_0}a_0\over 1+R^{2^{k+1}n_0}};$$
summing these inequalities for all odd $n_0$ we get the desired conclusion.

\hfill$\blacksquare$

\bigskip

\begin{prop} 1) Let $R_0$ the unique solution in $]0,1[$ of the following equation :
$$\sum_{n=1}^{\infty}R^{n}\vert a_n\vert\leq \sum_{n=0[2],n\geq 1}{4R^{n}\over
  1+R^{2n}}+\sum_{n=1[2]}{4R^{n}\over 1-R^{2n}}=1.$$ 
  Then, we will  have Bohr's phenomenon if  $R\leq R_0$,
for all $f\in\mathscr O(\mathscr E,\mathbb D)$.

\medskip
\noindent 2) Let $R_1$ be the unique solution in $[0,1]$ of
 $$\sum_{1}^{\infty}{4R^n\over
(1+R^{2n})}=1$$ 
  Then, we will  have Bohr's phenomenon if  $R\leq R_1$,
for all  holomorphic functions $f\in\mathscr O(\mathscr E,\mathbb D)$.
  with reals coefficients.
\end{prop}

\bigskip

\noindent\textbf{Proof : } 1)  Let $f=a_0+\sum_{1}^{\infty}a_n
F_{n,\mathscr E}\in\mathscr O(\mathscr E,\mathbb D)$. Up to a rotation we have  $a_0\geq 0$. Consider
$g=1-f$, she satisfies $\texttt{re}(g)>0$ and we can applies to $g$ all the preceedings results. We will have  Bohr's phenomena if we can find $R\leq r\leq 1$ such that  
$$\vert a_0\vert +\sum_{1}^{\infty}\vert a_n\vert \left(r^{n}+{R^{2n}\over r^n}\right)\leq
1,$$ %(remember that the "sup-norm" of  $\vert F_n\circ\Phi^{-1}\vert$ on $\{\vert z\vert =r\}$
%is $(r^{n}+{1\over r^n})$). 
The left side of this inequality is an increasing function of $r$, so such an inequality will be possible if 
$$\vert a_0\vert +\sum_{1}^{\infty}2\vert a_n\vert R^{n}\leq
1.$$ But, because the lemma 4.4 :
$$\vert a_0\vert +\sum_{1}^{\infty}2\vert a_n\vert R^{n}\leq \vert a_0\vert+\sum_{n=0[2],n\geq 1}{4R^{n}(1-a_0)\over
  1+R^{2n}}+\sum_{n=1[2]}{4R^{n}(1-a_0)\over 1-R^{2n}},$$ 
  and so, if
$$\sum_{n=0[2],n\geq 1}{4R^{n}\over
  1+R^{2n}}+\sum_{n=1[2]}{4R^{n}\over 1-R^{2n}}\leq 1,$$ 
 we will  assure the existence of Bohr's ph\'enomena. 
 
 \medskip
 \noindent 2) If the coefficients $a_n$ are reals we then can use the inequality $\vert a_n\vert\leq \frac{2 \texttt{re}(a_0)}{1+R^{2n}},\ n\geq 1$ (observed in (footnote 5)  the proof of the lemma 3.2). The result follow  immediatly.  \hfill $\blacksquare$

 \vfill\eject
\section{Optimality}

\subsection{Strategy}

\Large 
In this paragraph, we construct families of holomorphic functions $\phi_1(r,z)$ et 
$\phi_2(r,z)$ which gives optimality for the Bohr radius of the ellipse in the category of holomorphic functions  with arbitrary coefficients and also\footnote{Note that theses two radius are equal for the disc.} in the category of holomorphic functions with real coefficients.
One more time let $F_{n,\mathscr E}$ be the Faber polynomials of the ellipse\footnote{Remember that
  we have $F_{n,\mathscr E}(\Phi^ {-1}(z))=z^n+{R^{2n}\over z^n}$ et $F_{0,\mathscr E}=1$}, this is an orthogonal (not orthonormal) family of polynomials 
for the image measure on the boundary of the ellipse of the Lesbesgue measure on the unit circle via $\phi^{-1}$. Let us now consider the Bergman function associated :
$\sum_n\overline{F_{n,\mathscr E}(w_0)}F_{n,\mathscr E}(z)$ where $w_0\in\partial\mathscr E$ is fixed. 

To define extremal functions for Bohr's problem on the ellipse, the idea is to take sequences of points $(w_{0}^k)_k$ inside  the ellipse which tends to the boundary point $w_0$ (observe that this is the same in classical cases of the unit disc) and to perturb the family of Bergman function associated $\sum_n\overline{F_{n,\mathscr E}(w_0^ k)}F_{n,\mathscr E}(z)$. Because of the geometry of the ellipse, it seems reasonable to expect that we should choose the boundary points $w_0$ also on the axes of the ellipse and chossing  the sequences $(w_{0}^k)_k$ associated tending on the semi-axes to the boundary points. And that's really whats occurs as we soon shall see. 

 Clearly, such an asymetry doesn't accurs for the disc. Observe also that in the cases of the disc (i.e. $R=0$) we fall down on the classical functions giving optimality.
Thats what we gets when choosing the sequences $(w_{0}^k)_k$

\subsection{Somme technical lemmas}

Fix $0<R<1$ and consider for $R<r<1$ the function 
$$\phi_1(r,z)=-r+{1+r\over \gamma(r)}\sum_{n=1}^{\infty}{r^n+{R^{2n}  r^{-n}}\over
(1+R^{2n})^2}(z^n+{R^{2n} z^{-n}}),$$ 
where
$$\gamma(r)=\sum_{n=1}^{\infty}{r^n+{R^{2n} r^{-n}}\over 1+R^{2n}},$$

\bigskip

Let $(r_k)_k$ a real sequence converging to $1$ and consider the complex sequence $(z_k)_k$ defined by
 $$\vert\phi_1(r_k, z_k)\vert:=\sup_{\vert z\vert =1}\vert
\phi_1(r_k,z)\vert.$$ Up to to replace  $(r_k)_k$ by a  a subsequence, we can suppose $(z_k)_k$ converge, say to $z_0\in\partial\mathbb D$. 

\bigskip
In a same spirit, define for all $R<r<1$ :
$$\begin{aligned}\phi_2(r,z)&=-r+{1+r\over
\theta(r)}\biggl(\sum_{n=0[2],n\geq 1}{i^n(r^n+{R^{2n} 
r^{-n}})(z^n+{R^{2n}  z^{-n}})\over (1+R^{2n})^2}\\
&\hspace{3cm}-\sum_{n=1[2],n\geq 1}{i^n(r^n-{R^{2n} r^{-n}})(z^n+{R^{2n} 
z^{-n}})\over (1-R^{2n})^2}\biggl),
\end{aligned}$$ 
where
$$\theta(r)= \sum_{n=0[2],n\geq 1}{(r^n+{R^{2n} 
r^{-n}})\over (1+R^{2n})}+\sum_{n=1[2],n\geq 1}{(r^n-{R^{2n} 
r^{-n}})\over (1-R^{2n})} ,$$

\medskip
The sequence $(z_k)_k$ being defined as for the  $(\phi_1)_r$. We have :

\bigskip
\begin{prop}
$$\lim_{k\to\infty}{\vert\phi_1(r_k, z_k)\vert^{2}-1\over
1-r_k}=0,\quad {\text{and}}\quad \lim_{k\to\infty}{\vert\phi_2(r_k, z_k)\vert^{2}-1\over
1-r_k}= 0.$$ 
\end{prop} 

\bigskip
This clearly implies that 
$$\lim_{k\to\infty}{\vert\phi_1(r_k, z_k)\vert-1\over
1-r_k}= 0,\quad {\text{and}}\quad\lim_{k\to\infty}{\vert\phi_2(r_k, z_k)\vert-1\over
1-r_k}= 0.$$

\bigskip
We will prove the proposition 5.1 in the next paragraph. Before, we need some technical lemmas.

\bigskip
\begin{lem} We have the following estimations :
$$\begin{aligned}(1-r)\gamma(r)&=(1-r)\sum r^n+(1-r)\epsilon_1(r)=r+(1-r)\epsilon_1(r),\\
(1-r)\theta(r)&=(1-r)\sum r^n+(1-r)\epsilon_2(r)=r+(1-r)\epsilon_2(r),
\end{aligned}$$
where $\lim_{r\to 1}\epsilon_1(r)=\lim_{r\to 1}\epsilon_2(r)=0$. 
\end{lem}

\bigskip
\noindent\textbf{Proof : }\begin{enumerate}\item Straight computation gives
$$\gamma(r)-\sum_{n=1}^{\infty}r^n=\sum_{1}^{\infty}{{R^{2n}  r^{-n}}-r^nR^{2n}\over 1+R^{2n}},$$
and the left side of the equality is real analytic on a neighborought of  $r=1$ because $R<1$ and takes value $0$ if $r=1$. The result follows.

\item  Similarly  :
$$\theta(r)-\sum_{n=1}^{\infty}r^n=\sum_{n=0[2]}{{R^{2n}  r^{-n}}-r^nR^{2n}\over 1+R^{2n}}
+\sum_{n=1[2]}{r^nR^{2n}-{R^{2n} r^{-n}}\over 1-R^{2n}},$$
and as in the first cases, the right part of the equality  is real analytic on a neighborought of  $r=1$ because $R<1$ and takes value $0$ if $r=1$ ; this gives the result. \hfill$\blacksquare$
\end{enumerate}

\bigskip
For all $k\geq 1$, let us fix the following notations :
$$\begin{aligned}
A_k+iB_k&=\sum_{n=1}^{\infty}{r_{k}^n+{R^{2n} r_{k}^{-n}}\over
(1+R^{2n})^2}(z_{k}^n+{R^{2n}  z_{k}^{-n}}),\\
 C_k+iD_k&=\sum_{\underset{n=0[2],}{n\geq 1}}{i^n(r_{k}^n+{R^{2n} 
r_{k}^{-n}})(z_{k}^n+{R^{2n}  z_{k}^{-n}})\over (1+R^{2n})^2}\\
&\hspace{3cm}-\sum_{\underset{n=1[2],}{n\geq 1}}{i^n(r_{k}^n-{R^{2n}  r_{k}^{-n}})(z_{k}^n+{R^{2n} 
z_{k}^{-n}})\over (1-R^{2n})^2}.
\end{aligned}$$

\bigskip
\begin{lem} Write  :
$$\begin{aligned} 
A_k&=\texttt{re}\left(\sum_{n\geq 1} z_{k}^n r_{k}^n\right)+\alpha_k,\\
C_k&=\texttt{re}\left(\sum_{n\geq 1,n=0[2]}(ir_kz_k)^n-\sum_{n\geq
1,n=1[2]}(ir_kz_k)^n\right)+\beta_k.
\end{aligned}$$ 
Then  $\displaystyle \lim_{k\to\infty}\alpha_k=0=\lim_{k\to\infty}\beta_k$.

Moreover

-- If  $\lim_{k\to\infty}z_k=z_0\not =1$,then, there exists a constant $M_1>0$
such that  $\vert A_k+iB_k\vert\leq M_1$ for all  $k$ large enough. 

-- If $\lim_{k\to\infty}z_k=z_0\not =i$, then, there exists a constant $M_2>0$
such that  $\vert C_k+iD_k\vert\leq M_2$ for all  $k$ large enough.
\end{lem}

\bigskip
 
\noindent\textbf{Proof : } We have :
$$\begin{aligned}
A_k+iB_k&-\sum_{1}^{\infty}z_{k}^nr_{k}^{n}\\
&=\sum_{1}^{\infty}
{r_{k}^{n}{R^{2n}  z_{k}^{-n}}+z_{k}^{n}{R^{2n}  r_{k}^{-n}}+
{R^{4n} r_{k}^{-n}z_{k}^{-n}}-2R^{2n}r_{k}^{n}z_{k}^{n}-R^{4n}r_{k}^{n}z_{k}^{n}\over (1+R^{2n})^2}.
\end{aligned}$$
One more time, because  $R<1$,  the function on the right side   is real analytic on a neighborought of $r=1$ and $z=z_0$ with
$\vert z_0\vert=1$, so   is bounded for $k$ large enough. Function $\sum_{1}^{\infty}z_{k}^nr_{k}^{n}$ is also  ( for $k$ large enough) bounded if  $z_0\not =1$. These two observations assure the second part of the lemma for $A_k+iB_k$.

Moreover, observe that $\vert z_0\vert=1$ implies that the real part of the same function on the right side of the equality tends to $0$ as  $k\to+\infty$. This is the first part of the lemma for $A_k+iB_k$. Moreover, for $z_0=1$ the function itself tends to $0$ thats gives the lemma for
$A_k+iB_k$.

\medskip

We have the identity :
$$\begin{aligned} &C_k+iD_k-\sum_{n\geq 1,n=0[2]}(r_kiz_k)^n+\sum_{n\geq
1,n=1[2]}(r_{k}iz_{k})^n\\
&=\sum_{n=0[2]}{r_{k}^{n}{R^{2n}  (iz_{k})^{-n}}+(iz_k)^n{R^{2n} r_{k}^{-n}}
+{R^{4n} (r_kiz_k)^{-n}}-2R^{2n}(r_kiz_k)^n-R^{4n}(r_kiz_k)^n\over (1+R^{2n})^2}\\
&+\sum_{n=1[2]}{r_{k}^{n}{R^{2n}  (iz_{k})^{-n}}+(iz_k)^n{R^{2n}\over r_{k}^{-n}}
-{R^{4n}  (r_kiz_k)^{-n}}-2R^{2n}(r_kiz_k)^n+R^{4n}(r_kiz_k)^n\over (1-R^{2n})^2}.
\end{aligned}$$
he function on the right side   is real analytic on a neighborought of $r=1$ and $z=z_0$ with
$\vert z_0\vert=1$, (because $R<1$) so   is bounded for $k$ large enough.  
The function $\sum_{n\geq 1,n=0[2]}(r_kiz_k)^n-\sum_{n\geq 1,n=1[2]}(r_{k}iz_{k})^n$ is also bounded (for $k$ large enough) if $z_0\not =i$. These two observations implies the second part of the lemma for $C_k+iD_k$.

Note also that  $\vert z_0\vert=1$ implies that the real part of the same function on the right side of the equality tends to $0$ as  $k\to+\infty$. This is the first part of the lemma for $C_k+iD_k$. Moreover, if $z_0=i$, the function itself tends to $0$ thats gives the lemma for
$C_k+iD_k$. \hfill$\blacksquare$

 \bigskip

The two properties in the preceeding lemma means for $A_K+iB_k$ and $C_k+iD_k$ :

\bigskip

-- If $\lim_{k\to\infty}z_k=z_0=1$, then :
$$A_k+iB_k=\sum_{n\geq 1} r_{k}^nz_{k}^n+\lambda_k,$$
with $\lim_{k\to\infty}\lambda_k=0$.

-- If $\lim_{k\to\infty}z_k=z_0=i$, then 
$$C_k+iD_k=\sum_{n\geq 1,n=0[2]}(r_kiz_k)^n-\sum_{n\geq
1,n=1[2]}(r_{k}iz_{k})^n+\nu_k,$$ 
with $\lim_{k\to\infty}\nu_k=0$.

\bigskip
  \subsection{The proof of proposition 5.1}

\bigskip

Now, we can write :
$${\vert\phi_1(r_k, z_k)\vert^{2}-1\over
1-r_k}=-(1+r_k)-2{r_k(1+r_k)\over (1-r_k) \gamma(r_k)}A_k +
{(1+r_k)^2(1-r_k)\over (1-r_k)^2 \gamma^2(r_k)}(A_k^2+B_k^2),$$
$${\vert\phi_2(r_k, z_k)\vert^{2}-1\over
1-r_k}=-(1+r_k)-2{r_k(1+r_k)\over (1-r_k) \theta(r_k)}C_k +
{(1+r_k)^2(1-r_k)\over (1-r_k)^2 \theta^2(r_k)}(C_k^2+D_k^2),$$
So, to prove proposition 5.1, it is sufficient to show  
\begin{align}
\lim_{k\to\infty} -2{r_k(1+r_k)\over (1-r_k) \gamma(r_k)}A_k + {(1+r_k)^2(1-r_k)\over
(1-r_k)^2 \gamma^2(r_k)}(A_k^2+B_k^2)=2, \\
\lim_{k\to\infty}-2{r_k(1+r_k)\over (1-r_k) \theta(r_k)}C_k + {(1+r_k)^2(1-r_k)\over
(1-r_k)^2 \theta^2(r_k)}(C_k^2+D_k^2)=2.\end{align}

\bigskip

\noindent $\bullet$ First let us prove  (8) with $\lim_{k\to\infty}z_k=z_0\neq 1$.
 Because of lemma 5.3 :
$$\begin{aligned}\lim_{k\to\infty}&\left(-2{r_k(1+r_k)\over (1-r_k) \gamma(r_k)}A_k + {(1+r_k)^2(1-r_k)\over
(1-r_k)^2 \gamma^2(r_k)}(A_k^2+B_k^2)\right)\\
=&\lim_{k\to\infty}-2{r_k(1+r_k)\over (1-r_k) \gamma(r_k)}A_k=-4\lim_{k\to\infty}{A_k}.
\end{aligned}$$  
But we have
$$\lim_{k\to\infty}{A_k}=\lim_{k\to\infty}\texttt{re}\left(\sum_{n=1}^{\infty}r_{k}^{n}z_{k}^{n}\right)
=\lim_{k\to\infty}\texttt{re}\left({r_kz_k\over 1-r_kz_k}\right)={\texttt{re}(z_0)-1\over 2-2\texttt{re}(z_0)}=-{1/2},$$
what we had to prove.

\bigskip

\noindent $\bullet$ For   (9)  with $\lim_{k\to\infty}z_k=z_0\neq i$.

\noindent  Again because of lemma 5.3, it is sufficient to prove that
$$\lim_{k\to\infty}-2{r_k(1+r_k)\over (1-r_k) \theta(r_k)}C_k=\lim_{k\to\infty}-4{C_k}=2.$$
This is the case because (see below)
 $$\lim_k{C_k}=-\lim_{k\to\infty}\texttt{re}\left({r_kz_k\over r_kz_k-i}\right)=-{1/2},$$

\bigskip

\noindent $\bullet$ Now let us look at   (9)  with $\lim_{k\to\infty}z_k=z_0=i$.
\medskip

Let $c_k+id_k=C_k+iD_k-\nu_k$. Remember that :
$$C_k+iD_k=\sum_{n\geq 1,n=0[2]}(r_kiz_k)^n-\sum_{n\geq
1,n=1[2]}(r_{k}iz_{k})^n+\nu_k,$$ 
so
$$c_k+id_k:=\sum_{n\geq 1,n=0[2]}(ir_kz_k)^n-\sum_{n\geq
1,n=1[2]}(ir_{k} z_{k})^n.$$

After elementary computations we have the following equalitites  (where $w_k:=\texttt{im}{(z_k)}$) :
$$\begin{aligned} c_k&=-\texttt{re}\left({r_kz_k\over r_kz_k-i}\right),\\
\vert c_k+id_k\vert^2&={r_{k}^2\over 1-2r_kw_k+r_{k}^2},
\end{aligned}$$
so, the key equalitites :
$$ 
c_k=-{r_k(r_k-w_k)\over 1-2r_kw_k+r_{k}^2 },\quad
\vert c_k+id_k\vert^2={r_{k}^2\over 1-2r_kw_k+r_{k}^2}.
$$
Write in polar coordinates :
$$w_k=1+\rho_k\cos(\Lambda_k),\ r_k=1+\rho_k\sin(\Lambda_k),$$ 
where
$\rho_k\geq 0$ and $\Lambda_k\in [{\pi\over 2},{3\pi\over 2}]$ because $r_k\leq 1$ (the same for $w_k)$. We get the following

\begin{align}1-2r_kw_k+r_{k}^2=-2\rho_k\cos(\Lambda_k)-2\rho_{k}^2\cos(\Lambda_k)\sin(\Lambda_k)+\rho_{k}^2\sin^2(\Lambda_k),\\
2(r_k-w_k)+(1-r_{k}^2)(2-r_k)=-2\rho_k\cos(\Lambda_k)+\rho_{k}^2\sin^2(\Lambda_k)+o(\rho_{k}^2).\end{align}
From these, we can deduce that for $k$ large enough :
\begin{equation}1-2r_kw_k+r_{k}^2\geq a\rho_{k}^2,\end{equation} 
where $a>0$ is a constant. And
\begin{equation}-2\rho_{k}^2\cos(\Lambda_k)\sin(\Lambda_k)= \big(-2\rho_k\cos(\Lambda_k)+\rho_{k}^2\sin^2(\Lambda_k)\big)\mu_k,\end{equation} 
where $\lim_k\mu_k=0$. Using (12),  lemma 5.4 and 5.2,   we can replace in (9), $C_k$ and $D_k$ by $c_k$ and $d_k$.

Now, because of (9) and lemma 5.2, we have to prove that
$$\lim_{k\to\infty}\left({-2c_kr_k+{(1-R_{k}^2)\over (1-r_k)\theta(r_k)}\vert c_k+id_k\vert^2}\right)=1$$
or, using the lasts expressions for $c_k$ and $d_k$ and always lemma 5.2 :
$$\lim_{k\to\infty} {2(r_k-w_k)+(1-r_{k}^2)(1+(1-r_k)+(1-r_k)\epsilon_2(r_k))
\over 1-2r_kw_k+r_{k}^2} =1.$$
Because of (12), this limit is the same as
$$\lim_{k\to\infty}{2(r_k-w_k)+(1-r_{k}^2)(2-r_k)\over 1-2r_kw_k+r_{k}^2},$$
and because of (11) this last one is equal to
$$\lim_{k\to\infty}{-2\rho_k\cos(\Lambda_k)+\rho_{k}^2\sin^2(\Lambda_k)+o(\rho_{k}^2)\over 1-2r_kw_k+r_{k}^2}.$$
Then, we have the required conclusion by (12) and (13).

\bigskip

\noindent $\bullet$ Conclude with  (8) with $\lim_{k\to\infty}z_k=1$.
\medskip

Write $a_k+ib_k=\sum r_{k}^nz_{k}^n$. After elementary computations, we have the following (with $t_k=\texttt{re}{(z_k)}$) :
$$a_k=-{r_k(r_k-\texttt{re}{(z_k)})\over 1-2r_k\texttt{re}{(z_k)}+r_{k}^2 },\quad 
\vert a_k+ib_k\vert^2={r_{k}^2\over 1-2r_k\texttt{re}{(z_k)}+r_{k}^2},$$
which   assure that this case goes mutatis-mutandis as the last one, replacing $w_k$ by $\texttt{re}{(z_k)}$.\hfill$\blacksquare$

\subsection{Optimality : Functions with reals coefficients : }

\bigskip

Let us consider the family ${\phi_1(r_k,z)\over\phi_1(r_k,z_k)}$ on the unit disc. Their modulus less than $1$. Bohr's phenomenum on the ellipse will occurs only of there exists  $1>r_1>R$ such that
$${1\over\phi_1(r_k,z_k)}\left(r_k+{1+r_k\over \gamma(r_k)}\sum_{n=1}^{\infty}{r_{k}^n+{R^{2n}  r_{k}^{-n}}\over
(1+R^{2n})^2}\sup_{\vert z_1\vert=r_1}(z^n+{R^{2n}
z^{-n}})\right)\leq 1,$$ for all $k\in\mathbb N$. i.e.
$$(1+r_k)\sum_{n=1}^{\infty}{r_{k}^n+{R^{2n}  r_{k}^{-n}}\over
(1+R^{2n})^2}\sup_{\vert z_1\vert=r_1}(z^n+{R^{2n} z^{-n}})\leq
(\phi_1(r_k,z_k)-r_k)\gamma(r_k);$$
which leads to the existence of
$R<r_1\leq 1$, such that :
$$(1+r_k)\sum_{1}^{\infty}{r_{k}^n+{R^{2n}  r_{k}^{-n}}\over
(1+R^{2n})^2}(r_1^n+{R^{2n} r_1^{-n}})\leq
(\phi_1(r_k,z_k)-r_k)\gamma(r_k),$$ fol all $k\in\mathbb N$.
\bigskip

Because of proposition 1.1, $\phi_1(r_k,z_k)=(r_k-1)\varepsilon (r_k)+1 $
with $\lim_k\varepsilon(r_k)=0$. From this equality, the lemma 5.2, letting $k$ goes  to infinity in the last inequality leads to :
$$\sum_{1}^{\infty}{2\over
(1+R^{2n})}(r_1^n+{R^{2n}  r_1^{-n}})\leq 1.$$
This inequality is possible only if
$$\sum_{1}^{\infty}{4R^n\over
(1+R^{2n})}\leq 1,$$ 
because $R<r_1\leq 1$. This implies $R\leq R_1$ with $R_1\cong 0.205328678165046$.
\bigskip

\subsection{Optimality : The general case} We   follow steep by steep the \og real coefficients cases\fg{} replacing
$\phi_1$ by $\phi_2$. 

Let us consider  the family  $({\phi_2(r_k,z)\over\phi_2(r_k,z_k)})_k$ of holomorphic functions on the unit disc, their modulus is less than $1$, so Bohr's phenomenum on the ellipse will occurs only of there exists  $R<r_1\leq 1 $ such that  for all $k\in\mathbb N$.
\begin{multline}\notag
{1\over\phi_2(r_k,z_k)}\biggl[ r_k+{1+r_k\over \theta(r_k)}\sum_{n=0[2]}{r_{k}^n+{R^{2n}  r_{k}^{-n}}\over
(1+R^{2n})^2}\sup_{\vert z_1\vert=r_1}(z^n+{R^{2n}z^{-n}})     \\
+\frac{1+r_k}{\theta(r_k)}\sum_{n=1[2]} \frac{r_{k}^n-R^{2n}  r_{k}^{-n}}{(1-R^{2n})^2}
\sup_{\vert z_1\vert=r_1} (z^n+R^{2n}z^{-n})\biggr] \leq 1
\end{multline}

This leads to the existence 
(en exprimant le Sup) of
$R<r_1\leq 1$, such that
$$(1+r_k)\sum_{n=0[2]}{r_{k}^n+{R^{2n} r_{k}^{-n}}\over
(1+R^{2n})^2}(r_1^n+{R^{2n} r_{1}^{-n}})$$
$$+(1+r_k)\sum_{n=1[2]}{r_{k}^n-{R^{2n} r_{k}^{-n}}\over
(1-R^{2n})^2}(r_1^n+{R^{2n} r_{1}^{-n}})\leq
(\phi_2(r_k,z_k)-r_k)\theta(r_k),$$  for all $k\in\mathbb N$. 
\bigskip

Because  proposition 5.1, $\phi_2(r_k,z_k)=(r_k-1)\varepsilon (r_k)+1 $
with $\lim_k\varepsilon(r_k)=0$. One more time, this equality, the lemma 5.2, and letting $k$ goes  to infinity in the last inequality leads to :
$$\sum_{n=0[2]}{2\over
(1+R^{2n})}(r_1^n+{R^{2n} r_{1}^{-n}})+\sum_{n=1[2]}{2\over
(1-R^{2n})}(r_1^n+{R^{2n} r_{1}^{-n}})\leq 1.$$
But this last inequality is possible only if
$$\sum_{n=0[2]}{4R^n\over
(1+R^{2n})}+\sum_{n=1[2]}{4R^n\over
(1-R^{2n})}\leq 1,$$ because $R<r_1\leq 1$. This implies $R\leq
R_0$.

\bigskip

We have proved that

\bigskip
\begin{theo} Let $R_1$ be the unique solution in $[0,1]$ of
 $$\sum_{1}^{\infty}{4R^n\over
(1+R^{2n})}=1$$ and $R_0$ the unique solution in $[0,1]$ of
$$\sum_{n=0[2]}{4R^n\over
(1+R^{2n})}+\sum_{n=1[2]}{4R^n\over
(1-R^{2n})}=1.$$

1) If  $R>R_1$ then, there are no Bohr's phenomenon    for the ellipse in the category of holomorphic functions with reals coefficients.

\medskip

2) If  $R>R_0$ then, there are no Bohr's phenomenon  for the ellipse in the category of holomorphic functions.
\end{theo}

\end{document}